\documentclass[runningheads]{llncs}

\usepackage{amssymb}
\usepackage{latexsym}

\begin{document}

\pagestyle{headings}

\mainmatter

\title{A new notion of transitivity for groups and sets of permutations}

\titlerunning{TRANSITIVE SETS OF PERMUTATIONS}

\author{William J. Martin  \and Bruce E. Sagan}

\authorrunning{DRAFT \today}

\institute{Department of Mathematical Sciences \\
and Department of Computer Science \\
Worcester Polytechnic Institute \\
\email{martin@wpi.edu}
\and
Department of Mathematics \\
Michigan State University \\
\email{sagan@math.msu.edu}}

\maketitle

\def\u#1{{\underline{#1}}}
\def\Tr{{\mbox {\rm Tr}}}
\newcommand{\BMA}{{\mathbb A}}
\newcommand{\reals}{{\mathbb R}}
\newcommand{\A}{{\cal A}}
\newcommand{\B}{{\cal B}}
\newcommand{\C}{{\cal C}}
\newcommand{\E}{{\cal E}}
\newcommand{\F}{{\cal F}}
\newcommand{\I}{{\cal I}}
\newcommand{\J}{{\cal J}}
\newcommand{\Po}{{\cal P}}
\newcommand{\T}{{\cal T}}
\newcommand{\Sym}{{\mbox{\sf Sym}}}
\newcommand{\Sm}{\Sym(m)}
\newcommand{\sh}{{\mbox{\sf sh}}}
\newcommand{\height}{{\mbox{\sf d}}}
\newcommand{\colsp}{{\mbox{\sf colsp}}}
\newcommand{\al}{\alpha}
\newcommand{\la}{\lambda}
\newcommand{\Om}{\Omega}
\newcommand{\dom}{\unrhd}
\newcommand{\isdom}{\unlhd}
\newcommand{\orbi}{\mathop{\rm orb}\nolimits}
\newcommand{\ipr}[2]{\langle{#1},{#2}\rangle}


\begin{abstract}
Let $\Omega=\{1,2,\ldots,n\}$ where  $n \ge 2$. 
The {\em shape} of an ordered  set partition 
$P=(P_1,\ldots, P_k)$ of $\Omega$ is the integer
partition $\lambda=(\lambda_1,\ldots,\lambda_k)$
defined by $\lambda_i = |P_i|$. 
Let $G$ be a group  of permutations acting on $\Omega$. For a fixed
partition $\lambda$ of $n$, we say that $G$ is {\em $\lambda$-transitive}
if $G$ has only one orbit when acting on partitions $P$ of shape $\la$.
A corresponding definition can also be given when $G$ is just a set.
For example, if $\lambda=(n-t,1,\ldots,1)$, then a $\lambda$-transitive
group is the same as a $t$-transitive permutation group and if
$\lambda=(n-t,t)$, then we recover the $t$-homogeneous permutation
groups.

In this paper, we use the character theory of the symmetric group $S_n$ 
to establish some structural results regarding $\lambda$-transitive
groups and sets.  In particular, we are able to generalize a theorem
of Livingstone and Wagner~\cite{livwag} about $t$-homogeneous groups.
We survey the relevant examples coming from groups. While it is
known that a finite group of permutations can be at most 5-transitive
unless it contains the alternating group,
we show that it is possible to construct
a non-trivial $t$-transitive set of permutations for each positive integer 
$t$.   We also show how these ideas lead to a split basis for the
association scheme of the symmetric group.
\end{abstract}

\section{Introduction}
\label{Sec:intro}

Throughout, let $n \ge 2$ be an integer and let $\Omega = \{1,\ldots, n\}$.
Suppose that $G$ is a subgroup of the symmetric
group $S_n$.  Recall that $G$ is {\em $t$-transitive} if it is
transitive on ordered $t$-tuples of distinct elements of $\Om$.  It is
trivial to see that if $G$ is $t$-transitive, then it is
$(t-1)$-transitive.  Similarly, $G$ is {\em $t$-homogeneous} if it is
transitive on $t$-element subsets of $\Om$.  The following is one of
the main theorems in a paper of Livingstone and Wagner.
\begin{theorem}[\cite{livwag}]
\label{lwth}
Let $G\le S_n$ be $t$-homogeneous where $t\le n/2$.  Then $G$ is also
$(t-1)$-homogeneous.  \qed
\end{theorem}

Our goal is to define a more general notion of transitivity so that
both the $t$-transitive and $t$-homogeneous results will be special
cases.  To this end, let an {\em 
ordered set partition of $\Omega$} be an ordered tuple 
$$ P = (P_1, \ldots, P_k )$$
of pairwise disjoint non-empty subsets of $\Omega$ whose union is $\Omega$. 
These partitions are also called {\em tabloids} and play an important
role in the representation theory of $S_n$~\cite{jameskerber,sagan}.
We will only need to consider the case where $|P_i| \ge |P_{i+1}|$ for
all $1 \le i < k$.   On the other hand, an {\em (integer) partition
$\la$ of $n$}, written $|\la|=n$, is a sequence
$$ \lambda = (\lambda_1, \ldots, \lambda_k) $$
of positive integers whose sum is $n$ and satisfying the monotonicity
condition
$$ \lambda_1 \ge \lambda_2 \ge \cdots .$$
The $\la_i$ are called the {\em parts} of $\la$.
The {\em shape} of the set partition $P$ is the integer partition 
$$ |P| = (|P_1|, \ldots, |P_k|). $$

Let $\lambda$ be a fixed partition of $n$.
A subgroup $G$ of $S_n := {\sf Sym}(\Omega)$ is said to be {\em
$\lambda$-transitive} if it is transitive on partitions $P$ of shape
$\la$.   Equivalently, for any two set partitions $P$ and $Q$ 
with $|P|=|Q|=\lambda$, there exists $g \in G$ taking $P$ to $Q$ in
the sense that 
$$ g P_i = Q_i $$
(as a set) for each $i$, $1 \le i \le k$. This definition includes 
that of a $t$-transitive permutation group (use $\lambda=(n-t,1,\ldots,1)$)
and that of a $t$-homogeneous permutation group (use $\lambda=(n-t,t)$).

We extend this definition to arbitrary sets of permutations as
follows. For $D \subseteq S_n$ and $\lambda$ a
partition of $n$, we say that $D$ is {\em $\lambda$-transitive}
if there is a constant $r_\la$ such that, for any two set
partitions $P$ and $Q$ of shape $\la$, $D$ contains exactly
$r_\la$ permutations $g$ taking $P$ to $Q$.  We will also use the
terminology ``$t$-transitive,'' respectively ``$t$-homogeneous,'' for sets of
permutations which are $(n-t,1,\ldots,1)$-transitive, respectively
$(n-t,t)$-transitive.

The following two propositions follow directly from the definitions.
\begin{proposition}
If $G \le S_n$ is a $\lambda$-transitive group of permutations acting
on $\Omega_n$, then $G$ is also $\lambda$-transitive as a set. \qed
\end{proposition}

\begin{proposition}
If $D \subseteq S_n$ is a $\lambda$-transitive set of permutations
acting on $\Omega_n$ and $G=\langle D \rangle$ is the subgroup generated
by $D$, then $G$ is a $\lambda$-transitive group. \qed
\end{proposition}

Each set partition $P$ determines a Young subgroup $Y=Y_P$ which
consists of the set of permutations of $\Omega$ which fix each
subset $P_i$ of $P$ setwise. We say $Y$ has shape $\lambda$ where
$\lambda=|P|$.  There is a bijection between all set partitions of
shape $\la$ and all left cosets of $Y_P$ that preserves the action of
$S_n$.  Using this observation, it is easy to prove the following
characterization of $\la$-transitivity. 

\begin{lemma} \label{Ltransiffdesign}
Suppose $|\lambda|= n$ and $D \subseteq S_n$. Then $D$ is 
$\lambda$-transitive if and only if there is a constant 
$r_\la$ such that for any coset $gY_P$ of
any Young subgroup of shape $\lambda$ we have $|gY_P \cap D| = r_\la$. \qed
\end{lemma}

The last definition that will be necessary to understand our
generalization of the Livingstone-Wagner Theorem and 
($t$-transitivity result) 
is that of the
dominance order on integer partitions.
If $|\la|=|\mu|=n$ then
{\em $\lambda$ is dominated by $\mu$},
written $\lambda \unlhd \mu$, if and only if 
$$ \sum_{i=1}^j \lambda_i \le \sum_{i=1}^j \mu_i \qquad 
{\mbox {\rm (for all $j$)}} $$ 
where we define $\lambda_i=0$ if $i$ is greater than the number of 
parts in $\lambda$ and similarly for $\mu$.  We will give two proofs
that if $\mu\dom\la$ then $\la$-transitivity implies
$\mu$-transitivity, from which Theorem~\ref{lwth} follows immediately.
The first, which is just for groups but also includes a stronger
result about orbits, is essentially a generalization of the one given
by Livingstone and Wagner themselves.  The second proof holds for all
$\la$-transitive sets $D$.  But if $D$ is not a group then it does not
have orbits and so a combinatorial approach is needed in this case.

\section{The group case}
\label{Sec:group}

Let the permutation group $G$ act on two sets $S$ and $T$ with
characters $\chi$ and $\psi$, respectively.  The inequality
$\chi\ge\psi$ means that $\chi=\psi+\phi$ where $\phi$ is  some
character of $G$.  Equivalently, if $\chi^1,\ldots,\chi^k$ are the
irreducible characters of $G$ and we write $\chi=\sum_i c_i\chi^i,
\psi=\sum_i d_i\chi^i$ for nonnegative integers $c_i,d_i$ then
$\chi\ge\psi$ if and only if $c_i\ge d_i$ for all $i$.  We will also
need the notation
$$
\orbi(G,S)=\mbox{the number of orbits of $G$ acting on $S$.}
$$

\begin{lemma}
\label{Lorb}
Let $G$ be a permutation group acting on sets $S$ and $T$ with
characters $\chi$ and $\psi$, respectively.  If $\chi\ge\psi$ then
$$
\orbi(G,S)\ge\orbi(G,T).
$$
\end{lemma}
\proof
Taking inner products with the trivial character and using
$\chi=\psi+\phi$, we have
$$
\orbi(G,S)=\ipr{\chi}{1}=\ipr{\psi+\phi}{1}=\ipr{\psi}{1}+\ipr{\phi}{1}
\ge\ipr{\psi}{1}=\orbi(G,T).\qed
$$

We now apply this lemma to the situation at hand. 
Let $P$ be any set partition of shape $\la$ and let $S_\la$ denote the
set of all left cosets of the Young subgroup $Y_P$.  Of course,
$S_\la$ depends on $P$ and not just $\la$, but that fact will not play
a role in the following theorem.
\begin{theorem}
\label{Torblamu}
If $G\le S_n$ then
$$
\la\isdom\mu \Longrightarrow \orbi(G,S_\la)\ge\orbi(G,S_\mu).
$$
\end{theorem}
\proof
Let $\chi$ and $\psi$ be the characters of $G$ acting on $S_\la$ and
$S_\mu$, respectively.  Then by the previous lemma, it suffices to show
that $\chi\ge\mu$.  Let $P$ and $Q$ be ordered set partitions giving
rise to the sets $S_\la$ and $S_\mu$, respectively.  Then
$\chi$ and $\psi$ are the restrictions to $G$ of the characters
$\chi'=1_{S_\la}\uparrow^{S_n}$ and $\psi'=1_{S_\mu}\uparrow^{S_n}$,
respectively, so we are reduced to showing $\chi'\ge\psi'$.  But
then by Young's Rule
$$
\chi'=\sum_{\nu\dom\la} K_{\nu\la}\chi^{\nu}
\quad\mbox{and}\quad
\psi'=\sum_{\nu\dom\mu} K_{\nu\mu}\chi^{\nu}
$$
where $\chi^\nu$ is the irreducible character of $S_n$ corresponding
to $\nu$ and $K_{\nu\la}$ is a Kostka number.  However, it is
known~\cite{lv:opc,white}
that if $\la\isdom\mu$ then $K_{\nu\la}\ge K_{\nu\mu}$ for all $\nu$.
This finishes the proof. \qed

As an immediate corollary, we get our first generalization of the
Livingstone-Wagner Theorem.
\begin{theorem}
\label{TGlamu}
If $G\le S_n$ is $\la$-transitive and $\mu\dom\la$ then $G$ is
$\mu$-transitive. \qed
\end{theorem}

\section{The center of the group algebra and the inner distribution}
\label{Sec:center}

Our main tool for extracting information about $\la$-transitive sets
will be the center of the group algebra of the 
symmetric group $S_n$, specifically in the form arising from
the left regular representation  $g \mapsto A(g)$ 
(over the reals).   We will denote this commutative algebra of
symmetric matrices by $\BMA$ and 
its dimension by $p(n)$, the number of partitions of the integer $n$. 
One may learn quite a bit about this algebra via the character theory of the 
symmetric group; see references~\cite{banito,jameskerber,sagan}.
We identify two distinguished bases of $\BMA$ as follows.
The first basis contains one matrix for each 
conjugacy class $\C_{\alpha}$, namely
$$ A_{\alpha} = \sum_{g \in \C_{\alpha} } A(g). $$
The second distinguished basis is the  basis of primitive idempotents for the center.
This basis contains one element, $E_{\mu}$ say, for each irreducible representation
of $S_n$. The two bases are related by the equations
$$ E_{\mu} = {f_{\mu}  \over n!} \sum_{\alpha} \chi_{\alpha}^{\mu} \, A_{\alpha} $$
where $f_{\mu}$ is the degree of the 
irreducible representation $\rho_{\mu}$ and $\chi_{\alpha}^{\mu}$ is the value
taken by the corresponding irreducible character on the 
conjugacy class $\C_{\alpha}$ \cite[Theorem~II.7.2]{banito}.

The matrices in $\BMA$ are clearly simultaneously diagonalizable.
The {\em eigenspace of $\BMA$ indexed by $\mu$} is the column
space of $E_{\mu}$:
$$ V_{\mu} = \colsp (E_{\mu}). $$
The dimension of $V_{\mu}$ is $f_{\mu}^2$.

\bigskip

Let $D$ be a non-empty subset of $S_n$ and let $x = x_D$ 
denote the indicator
vector of $D$. This is a 01-vector of length $n!$ with $x(g) = 1$ if and only
if $g  \in D$. The {\em inner distribution} of  $D$ is the rational vector of
length $p(n)$ with $\alpha$ entry equal to 
$$ a_{\alpha} = {1 \over |D|} x^\top A_{\alpha} x. $$
The  {\em dual distribution} of $D$ 
is the vector  $b$ with entries indexed by irreducible representations and 
defined by 
\begin{equation} \label{EEmubmu}
b_{\mu} = {n! \over |D|} x^\top E_{\mu} x . 
\end{equation}
Since each $E_{\mu}$ is positive semi-definite, we have $b \ge 0$ and $b_{\mu}=0$ if and
only if $E_{\mu} x = 0$. 
And since $E_{\mu}$ is the 
matrix representing orthogonal projection
onto $V_{\mu}$ we have another equivalent condition, namely
$b_{\mu}=0$  if and only if $x \bot V_{\mu}$.

For arbitrary $D$, the quantities $a_{\alpha}$ and $b_{\mu}$ have 
alternative expressions as follows:
$$ 
a_{\alpha} = {1 \over |D|} \left|\{ (g,h) \in D \times D : gh^{-1} \in \C_{\alpha} \}\right|; 
$$
and
\begin{equation} \label{Ebmu}
b_{\mu} = {f_\mu \over |D|} \sum_{g, h \in D}\chi^{\mu}(gh^{-1}).
\end{equation}
When $D$ is a group, these simplify to
$$ a_{\alpha} = \left| D \cap \C_{\alpha} \right| $$
and 
\begin{equation} \label{Ebmugp}
b_{\mu} = f_{\mu} \sum_{g \in D } \chi^{\mu}(g). 
\end{equation}

\section{The set case}
\label{Sec:setcase}

\begin{lemma}
Suppose $D$ is a coset of a Young subgroup of shape $\lambda$ and $b$ is the
dual distribution of $D$. Then $b_{\mu} \neq  0$ if and only if
$\mu \unrhd \la$. 
\end{lemma}

\proof Because of equation~(\ref{Ebmu}), we may assume $D$ is a Young 
subgroup.  In this case, by Young's rule and~(\ref{Ebmugp}), 
the number $b_{\mu}$ is a non-zero multiple of the Kostka 
number $K_{\mu\lambda}$.   \qed

Let $Y_{\lambda}$ denote the incidence matrix of $S_n$
versus cosets of Young subgroups
of shape $\lambda$. The $n!$ rows of $Y_{\lambda}$ are indexed
by permutations and the columns are indexed by all left cosets $hY_P$ of all 
Young subgroups of shape $\lambda$.  In the entry indexed $(g,hY_P)$,
we have a $1$ if $g \in hY_P$ and a $0$ otherwise.

\begin{corollary} \label{CSlambda}
Let $W_\la$ be the column space of $Y_\la$.  Then
$$ W_\lambda= \bigoplus_{\mu \unrhd \lambda} V_{\mu}.$$
\end{corollary}

\proof  This space is invariant under the left action of the group
algebra.  So, in particular, it is $\BMA$ invariant and thus a sum of
eigenspaces.  Combining this with the previous lemma yields the
desired  result. 
\qed

\begin{theorem} \label{Tmain}
Let $D$ be a set of permutations in $S_n$ and $\lambda$ a partition of $n$. 
Then $D$ is $\lambda$-transitive if and only if for every $\mu$ with
$\lambda \unlhd \mu \lhd (n)$ we have $b_{\mu}=0$.
\end{theorem}

\proof  We prove the forward implication noting that the converse follows
by reversing the steps of the proof.

>From Corollary \ref{CSlambda}, the eigenspaces which appear in $W_\lambda$
are precisely those $V_\mu$ for which $\mu \unrhd \lambda$.
Now if $x$ denotes the indicator vector of $D$, then 
by Lemma~\ref{Ltransiffdesign} we have 
$x^\top Y_{\lambda} = r_{\lambda} {\bf 1}$ where ${\bf 1}$ stands for the 
all ones vector of appropriate length.  Let $x_0$ denote the projection of 
$x$ onto $W_{\lambda}$.  
Note that $x$, and thus $x_0$, has the same
inner products with columns of $Y_\la$ as does the vector
${r_{\lambda} \over |Y_P|}{\bf 1}$  where $P$ is of shape $\la$.
So the same will be true for products with the vectors from
any orthonormal basis $B$ for $W_{\lambda}$.
It follows that $x_0$ is equal to a multiple of the
all ones vector. Since $x_0\in V_{(n)}$, we have $x \bot V_{\mu}$ for
all $\mu$ such that $\lambda \unlhd \mu \lhd (n)$ and hence the
theorem.  \qed

An immediate corollary of this result is our second generalization of
the Livingstone-Wagner Theorem.

\begin{theorem} \label{TDlamu}
If the set $D$ is $\lambda$-transitive and $\mu \dom \lambda$, then $D$ is 
$\mu$-transitive. \qed
\end{theorem}

We get a second corollary of Theorem~\ref{Tmain} by combining it with
equation~(\ref{Ebmu}).  
\begin{theorem}
\label{Tchi}
The set $D$  is $\lambda$-transitive if and only if
$$ \sum_{g,h \in G} \chi^{\mu}(gh^{-1}) = 0 $$
for each $\mu$ such that $\lambda \unlhd \mu \lhd (n)$.  \qed
\end{theorem}

\section{Examples}
\label{Sec:eg}

>From Theorem~\ref{TGlamu}, we know that a $\lambda$-transitive group of 
permutations is $t$-homogeneous where $t=\sum_{i=2}^k \lambda_i$. Thus
we are severely restricted in our search for $\lambda$-transitive groups.
Other than alternating and symmetric groups, there are no $t$-homogeneous
groups for $t \ge 6$. (This is a consequence of the 
classification of finite simple groups; see \cite{dm} or \cite{handbk}.)  
There are precisely two
non-trivial $5$-transitive groups: the Mathieu groups $M_{12}$ and $M_{24}$.
Of course, a $t$-transitive group is 
$\lambda$-transitive for any partition $\lambda$ refined by 
$(n-t,1,\ldots,1)$. So we seek $t$-homogeneous
groups which fail to be $t$-transitive. Livingstone and Wagner~\cite{livwag}
showed that there are no such groups for $t\ge 5$. The balance of 
the following result is due to Kantor~\cite{kantor}  who dealt with the 
cases $t=2,3,4$.

\begin{theorem}[Kantor, Livingstone/Wagner \cite{dm}]
Suppose $G$ is a permutation group acting $t$-homogeneously on $n$ points ($t \le n/2$), 
but not $t$-transitively. Then one of the following is the case:
\begin{description}
\item[$\bullet$]  
$t=2$, $ASL_1(q) \le G \le A\Sigma L_1(q)$, $q \equiv 3 \bmod{4}$;
\item[$\bullet$]    
$t=3$,  $PSL_2(q) \le G \le P\Sigma L_2(q)$, $q \equiv 3 \bmod{4}$;
\item[$\bullet$]    
$t=3$, $G = AGL_1(8)$, $A\Gamma L_1(8)$, $A\Gamma L_1(32)$;
\item[$\bullet$]    
$t=4$, $G = PGL_2(8)$, $P\Gamma L_2(8)$, $P\Gamma L_2(32)$.
\end{description}
In all cases, the group acts $(t-1)$-transitively.
\end{theorem}

For each group appearing in this theorem, we wish to determine its
full level of transitivity.
(For $t=2$, there is no partition between $(n-2,2)$ and $(n-2,1,1)$ 
to consider.)
\begin{itemize}
\item For $PSL_2(q)$ acting on the projective line ($q+1$ points), 
we find the stabilizer of
a $3$-set  has size three. Since the action is faithful, we conclude
that $PSL_2(q)$ acts $(n-3,2,1)$-transitively on the projective line.
It follows that this is also the full degree of transitivity of 
$P\Sigma L_2(q)$ and any group between the two. 
\item  $AGL_1(8)$ has order 56, so this group acts sharply 
3-homogeneously on ${\mathbb F}_8$.
\item $  A\Gamma L_1(8)$ acts faithfully on ${\mathbb F}_8$ with
three group elements stabilizing any 3-set. So this action is
$(5,2,1)$-transitive.
\item $  A\Gamma L_1(32)$  has order 4960, so this group acts 
sharply 3-homogeneously on ${\mathbb F}_{32}$.
\item Finally, we consider the groups with $t=4$. First,
$PGL_2(8)$ has 504 elements and acts faithfully on $PG(1, 8)$ with
the stabilizer of any 4-set being isomorphic to the Klein four-group, so
it is $(5,3,1)$-transitive in this action.
\item $P\Gamma L_2(8)$ has 1512 elements and acts 
4-homogeneously on nine points. 
The stabilizer of a 4-set has size 12, and thus is isomorphic to $A_4$. 
Hence this group is $(5,2,1,1)$-transitive. 
\item For $P\Gamma L_2(32)$,
the stabilizer of a 4-set is isomorphic to the Klein four-group, so 
the action is $(29,3,1)$-transitive.
\end{itemize}

We note that, in 1975, Peter Neumann~\cite{neumann}  gave a 
classification scheme for $t$-homogeneous groups which are not 
fully $t$-transitive in terms of the isomorphism type 
of the stabilizer of a $t$-set. In Neumann's language, 
a group which acts  {\em almost generously $3$-transitively} 
is one in which the stabilizer of any $4$-set contains $A_4$. (He does
not presume  
a $4$-homogeneous action in this definition.) 
Clearly a group acting $4$-homogeneously has this property exactly when
it is $(n-4,2,1,1)$-transitive.  Neumann calls a  permutation group 
{\em a little generously $3$-transitive} if the stabilizer of any 
$4$-set contains the 
Klein four-group. 
Observe that a group acting $4$-homogeneously and having this property 
will be  $(n-4,3,1)$-transitive. Since there are no  $4$-homogeneous
groups with cyclic stabilizer of a $4$-set, the converse is true as
well.  But this observation hinges on Kantor's result.

\bigskip 

We now construct a family of $\lambda$-transitive sets of permutations 
for partitions $\lambda$ having the property that any $\lambda$-transitive
group must be strictly larger than the set we construct.

\begin{example} \cite{bieredel}
For $q$ an odd prime power, the affine group $AGL_1(q)$ acts 
2-transitively on ${\mathbb F}_q$.
There can be no $2$-homogeneous group of size ${q \choose 2}$ but there
is such a set of permutations. Choose a subset $S$ of 
${\mathbb F}_q - \{0\}$ such that $S$ is disjoint from $-S$ and 
$S \cup (-S) = {\mathbb F}_q - \{0\}$.
(Since $q$ is odd, this is always possible.) Let
$$ D = \{ x \mapsto ax + b :  a \in S, b \in {\mathbb F}_q \}. $$
We claim that this set is $(q-2,2)$-transitive.  For if the mapping 
$x \mapsto ax+b$ takes $u$ to $y$ and $v$ to $z$, then the mapping which
takes $u$ to $z$ and $v$ to $y$ (in $AGL_1(q)$) is $x \mapsto (-a)x+b'$ for
some $b'$. Thus no two elements of $D$ take an unordered pair $\{u,v\}$ to
the same image.  \qed
\end{example}

\section{A recursive construction}
\label{Sec:construct}

In this section, we show that, given any partition $\lambda$ of a 
positive integer, $m$ say,
if we are permitted to increase the largest part (thus
obtaining a partition $\lambda'$ of some integer $n \ge m$),
we can construct a $\lambda'$-transitive set of permutations $D$
which is ``small'' in comparison to $|S_n|$.

Suppose the following three structures are given:
\begin{description}
\item[$\bullet$] 
a $t$-$(n,k,\nu)$ block design on point set $\Omega = \{1,\ldots,n\}$
with block set $B$;
\item[$\bullet$]  
a $t$-transitive set of permutations $D_1 \subseteq S_k$;
\item[$\bullet$]  
a $t$-transitive set of permutations $D_2 \subseteq S_{n-k}$.
\end{description}
>From these, we will construct a $t$-transitive set $D$ of permutations
inside $S_n$.  

We first define a number of bijections. Let 
$$\tau : \{ k+1,\ldots,n\} \longrightarrow \{1,\ldots,n-k \} $$
via $j \mapsto j-k$. For each block $b \in B$, choose bijections
$$ \phi_b : \{ 1,\ldots, k\}  \longrightarrow b$$
and
$$ \psi_b :\{ 1,\ldots, n-k\} \longrightarrow \Omega - b , $$
which will remain fixed.

\begin{lemma}
With the above notation, the set of functions
$$ D = \left\{ (\phi_b \circ \pi) \cup (\psi_b \circ \sigma \circ \tau) :
b \in B, \pi \in D_1, \sigma \in D_2  \right\} $$
is a $t$-transitive set of permutations in $S_n$.
\end{lemma}

\begin{example}
Suppose $t=2$, $n=7$, $k=3$ and the blocks are those of the Fano plane
$$ B = \left\{ 124, 235, 346, 457, 561, 672, 713\right\}, $$
$D_1=S_3$ and $D_2=A_4$. We then obtain a 2-transitive set of $7.6.12=504$
permutations in $S_7$. For instance, 
if $b=\{3,4,6\}$ and $\Omega-b=\{5,7,1,2\}$
(listed in the order induced by $\phi_b$ and $\psi_b$) and $\pi =
(12)$, $\sigma=(13)(24)$, the resulting  
permutation in $D$   is 
$$ (\phi_b \circ \pi) \cup (\psi_b \circ \sigma \circ \tau) = (14)(2365).
 \qquad \qed $$
\end{example}

The following notation will be required in the proof of this 
lemma. Suppose $I$ and $O$ are disjoint subsets of $\Omega$.
We would like to count blocks $b \in B$ which contain all of the
points in $I$ yet none of the points in $O$. 
If $|I|=i$ and $|O|=j$ and $i+j\le t$, 
then there exists a
constant $\nu_{i,j}$ such that there exactly $\nu_{i,j}$ blocks
$b$ such that $I \subseteq b$ and 
$O \cap b =\emptyset$~\cite[Remark~4.17]{handbk}. 

Suppose that $h\le t$.  If $a_1,\ldots,a_h$ are distinct elements of
$\Omega_k$ and 
$b_1,\ldots,b_h$ are distinct elements of $\Omega_k$, then there
are exactly $c_h$ permutations $\pi \in D_1$ such that 
$\pi a_j = b_j$ for all $1\le j \le h$ for some constant $c_h$.
Similarly, there exist constants $d_0, \ldots, d_t$
associated to $D_2$.

{\em Proof of Lemma 5.}\
Let $a_1 < \cdots < a_t$ be distinct elements of $\Omega$
and let $b_1, \ldots, b_t$ be distinct elements of $\Omega$.
If $a_1 > k$, set $i =0$; otherwise, let $i$ be the largest 
integer such that $a_i \le k$ and set $j=t-i$.
There are $\nu_{i,j}$ blocks $b$ of the design $(\Omega,B)$
for which $b_1,\ldots,b_i \in b$ and $b_{i+1},\ldots,b_t \not\in b$.
There are $c_i$ permutations $\pi \in D_1$ taking $a_\ell$ to
$\phi_b^{-1}b_\ell$ 
for $\ell \le i$. There are $d_{j}$ permutations $\sigma$ in
$D_2$ taking $a_\ell - k$ to $\psi_b^{-1} b_\ell$ for $i < \ell \le t$.
Thus there are 
$$ r = c_i d_j \nu_{i,j} $$
permutations $\rho \in D$ satisfying $\rho a_\ell = b_\ell$ 
for all $1\le \ell \le t$.
To finish the proof, we simply must show that this quantity $r$  does not
depend on $i$. But we have
$$ \nu_{i+1,j} (n-i-j) = \nu_{i,j} (k-i) $$
and
$$ \nu_{i,j+1} (n-i-j) = \nu_{i,j} (n-k-j), $$
giving us
$$ \nu_{i+1,j-1} (n-k-j+1) = \nu_{i,j} (k-i). $$
Moreover, it is obvious that 
$c_i = c_{i+1} (k-i)$  and $ d_{j-1} = d_j (n-k-j+1)$.
Putting these identities together, we find that
$$ \nu_{i,j} c_i d_j = \nu_{i+1,j-1} c_{i+1} d_{j-1} $$
which completes our proof. \qed

We will now make use of the following result:

\begin{theorem}[Teirlinck~\cite{teir}]
Given integers $t$ and $n$ with 
$$n \equiv t \pmod{(t+1)!^{2t+1}}$$
and $n\ge t+1>0$, there exists a $t$-$(n,t+1,(t+1)!^{2t+1})$ 
design without repeated blocks. \qed
\end{theorem}

Using these designs in the previous construction, we obtain the following
asymptotic existence theorem.

\begin{theorem}
Let $\lambda_2 \ge \lambda_3 \ge \cdots \ge \lambda_t$ be a 
non-increasing sequence of positive integers and let
$\epsilon > 0$.  Then there exist infinitely many 
values  for  $\lambda_1 \ge \lambda_2$  such that there exists
a $\lambda$-transitive set $D$ of permutations inside $S_n$ 
($n=\lambda_1 + \cdots + \lambda_t$) having $|D|/n! < \epsilon$.
\end{theorem}

\proof It suffices to prove this for $\lambda_2 = \cdots = \lambda_t=1$.
Let $k=t+1$ and choose $n\ge 2t+1$ to satisfy Teirlinck's theorem. 
Then there exists a block design $B$ on $n$ points with each $t$-set
appearing in exactly $(t+1)!^{2t+1}$ blocks. 
Using $S_{t+1}$ and $S_{n-t-1}$ as our component $t$-transitive
permutation sets, we obtain
a $t$-transitive set $D$ of permutations in $S_n$ with 
$$ |D| = |B| \cdot (t+1)! \cdot (n-t-1)!  = { (t+1)!^{2t+1} \over t+1} \cdot 
{n \choose t} \cdot (t+1)! \cdot (n-t-1)!,  $$ 
that is, $|D| = K \cdot n! / (n-t)$ where $K = (t+1)!^{2t+1}$. 
With $t$ fixed, we see that 
$$ \lim_{n \rightarrow \infty} {|D| \over n!} = 
K \cdot \lim_{n \rightarrow \infty} {1 \over n-t}  = 0 . \qquad  \qed $$

\section{A split basis for the association scheme of the symmetric group}
\label{Sec:scheme}

We now turn our attention to the association scheme of the symmetric
group. The center $\BMA$ of the group algebra of $S_n$ is a vector space 
of symmetric matrices closed under both ordinary matrix multiplication
and under Schur (entrywise) multiplication and containing both $I$ and
the all-ones matrix $J$. In other words, $\BMA$ is a {\em Bose-Mesner
algebra}. The underlying association scheme has $S_n$ as its vertex set
and contains one graph $\Gamma_\alpha$ for each conjugacy class $\C_\alpha$:
two permutations $g$ and $h$ are adjacent in $\Gamma_\alpha$ if $gh^{-1}
\in \C_\alpha$. The adjacency matrix of $\Gamma_\alpha$ is the matrix $A_\alpha$
introduced in Section~\ref{Sec:center}. The basis of primitive central
idempotents $E_\mu$ was introduced earlier as well.

In the study of association schemes, one often chooses a base point from
the vertex set and considers the {\em trivial module} $\BMA {\bf e}$ 
where ${\bf e}$ is the elementary basis vector associated to the base point.
If we choose the identity permutation as our base point and use the
corresponding elementary basis vector ${\bf e}_{(1)}$,
then the trivial module $\BMA {\bf e}_{(1)}$ for $\BMA$
acting on $\reals^{n!}$ is simply the vector space of class functions.
The trivial module has two distinguished bases: the indicator functions
$A_{\alpha} {\bf e}_{(1)}$ for the conjugacy classes $\C_\alpha$ and 
the irreducible characters 
$\chi_\mu = \frac{n!}{f_\mu} E_\mu {\bf e}_{(1)}$.  We found earlier
that this module alone was sufficient for the study of $\lambda$-transitive
groups. Yet, to extend the results to arbitrary sets of permutations, we
found it necessary to appeal to the coding theorist's inner distribution and
dual distribution. This is, in effect, a smoothing technique in which
we average over all base points $g \in S_n$ in order to account for the
possible asymmetry of a set $D$.

Consider the set $\Po$ of all cosets of Young subgroups of 
$S_n$ partially ordered by
inclusion. We partition the elements of this poset into fibers corresponding
to the shape of the coset $gY_P$. The matrix $Y_\lambda$ introduced in
Section~\ref{Sec:setcase} is then the incidence matrix between 
the minimal elements (singletons) of this poset and the elements 
in the fiber indexed by $\lambda$. We are interested in the 
matrices $C_\lambda = Y_\lambda Y_\lambda^\top$ where 
$\lambda$ is a partition of $n$. 

\begin{theorem}
The set $\left\{ C_\lambda : |\lambda|=n \right\}$  is a basis for $\BMA$.
The change-of-basis equations with respect to the $A_\alpha$ have the form
$$
C_\lambda= \sum_{\alpha\isdom\la} m_{\lambda,\alpha}  A_\alpha 
$$
with $m_{\la,\la}>0$.  The change-of-basis equations with respect to
the $E_\mu$ have the form
$$ 
C_\lambda = \sum_{\mu\isdom\la} n_{\lambda,\mu} E_\mu
$$
with $n_{\la,\mu}\neq 0$ for all $\mu \isdom \lambda$.
\end{theorem}

\proof Since the entry of $C_\lambda$ in row $a$ column $b$ is 
equal to the number of cosets $hY_P$ of Young subgroups of shape $\lambda$
containing both $a$ and $b$,
we clearly have $A(g)^{-1} C_\lambda A(g) = C_\lambda$ for any $g \in S_n$
showing that $C_\lambda$ lies in $\BMA$. Next, note that the column
space of $C_\lambda$ is equal to that of $Y_\lambda$ so that,
by Corollary~\ref{CSlambda},
$n_{\lambda,\mu} \neq 0 $ if $\mu \dom \lambda$ and 
$n_{\lambda,\mu} =0$ otherwise.
So the matrices $C_\lambda$ indeed form a basis for $\BMA$. 

Now consider the coefficients $m_{\lambda,\al}$.  The $(a,b)$-entry 
of $C_\lambda$ is non-zero if and only if
there is some coset of some Young subgroup  of shape $\lambda$
which contains both $a$ and $b$. 
Equivalently, $ab^{-1}$ is an
element of that Young subgroup.
But the conjugacy classes $\C_\al$ which intersect a Young subgroup of 
shape $\la$ are precisely those for which $\al$ is a refinement of
$\la$.  So by definition of $A_\al$, $m_{\la,\al}\ge0$ with strict
inequality precisely when $\al$ refines $\la$.  The Theorem now
follows since refinement implies $\al\isdom\la$ and because $\la$
refines itself.  \qed

Following Terwilliger~\cite{itotanter,ter2,ter24}, 
we call the basis $\left\{ C_\lambda:
\lambda\right\}$ a {\em split basis} for the association scheme of $S_n$.
We remark, however, that the typical linear ordering on relations
and primitive idempotents is replaced in this case by the partial 
order $\unlhd$.

%
%
%

In the language of association schemes, a $\lambda$-transitive
set $D$ is the same as a Delsarte $\T$-design~\cite[p32]{delthesis} where 
$\T = \left\{ \mu : \mu \dom \lambda \right\}$.
This provides another viewpoint for our investigation and suggests
an avenue for finding future non-existence results. We remark
that Delsarte's linear programming bound --- the most obvious 
technique in this regard --- is mostly ineffective in this particular
setting,  being overshadowed by the following more elementary 
divisibility conditions on $\lambda$-transitive sets.

\begin{proposition} \label{Pdivis}
If $D \subseteq S_n$ is a $\lambda$-transitive set of permutations,
then $|D|$ is divisible by ${n \choose \mu}:={n \choose \mu_1,\ldots,\mu_k}$ 
for any partition $\mu \dom \lambda$.
\end{proposition}

\proof Since $D$ must also be $\mu$-transitive by Theorem~\ref{TDlamu},
each coset of each Young subgroup of shape $\mu$ contains $r_\mu$ permutations
from $D$ for some integer $r_\mu$.  
Given a fixed $Y_P$ of shape $\mu$, it has 
${n \choose \mu}$ cosets which partition $S_n$ and thus $D$.
So $|D| = r_\mu {n \choose \mu}$. 
\qed

\section{Open questions}
\label{Sec:open}

We list here some open questions raised by our work in the hopes that
the reader will be tempted to work on them.

(1) An interesting construction is given in \cite{conway} for a 
set of $13!/7!$ permutation in $S_{13}$ which Conway calls $M_{13}$.
While this set is six-transitive in the sense of \cite{conway},
it is not $(7,1,1,1,1,1,1)$-transitive  in the sense we are using
in this paper. We were unable to compute the full degree of 
transitivity of this set due to its size. We also ask if $M_{13}$ 
can be combined with some small collection of its translates to 
yield a $6$-transitive set of permutations in $S_{13}$.

(2)  Is there an analogue of $\la$-transitive groups and sets for
other Coxeter groups?  In particular, it would be interesting to
investigate the hyperoctahedral group.  One would hope that there is a
generalization of Theorem~\ref{TDlamu} in this setting.

(3) In view of Proposition~\ref{Pdivis}, 
no $\la$-transitive subset of $S_n$ can contain less than ${n \choose \la}$
permutations. Trivial examples where this lower bound is achieved are
the cyclic, alternating and symmetric groups. At the end of 
Section~\ref{Sec:eg}, we presented examples achieving this lower 
bound for $\la=(q-2,2)$ where $q$ is an odd prime power. We ask if 
there are any other infinite families of
examples achieving this lower bound.

(4) The {\em Krein parameters} $q_{\lambda,\mu}^\nu$ of the association
scheme of the symmetric group are defined by the equations
$$ E_\lambda \circ E_\mu  = \frac{1}{n!} \sum_\nu q_{\lambda,\mu}^\nu E_\nu $$
where $\circ$ denotes the entrywise product of matrices.
In the study of association schemes, one is often interested in
determining which Krein parameters vanish for a given family of schemes. 
While this question is known to be hard in general, we identify one
interesting result along these lines.
Using Theorem 2.9.22 on p99 of the book by James and 
Kerber~\cite{jameskerber}, we find
that the Krein parameter $q_{\lambda, \mu}^{\nu}$ vanishes if 
$\height(\nu) > \height(\lambda) + \height(\mu)$ where 
$\height(\nu) = \nu_2+\nu_3 + \cdots$.  

\section*{Acknowledgments}

The authors wish to thank Peter Cameron, Jon Hall, Robert Liebler,
Spyros Magliveras and Paul Terwilliger for helpful comments. 
This paper was written while the first author was visiting the Department 
of Combinatorics and Optimization at the University of Waterloo. He 
extends thanks to the department for its hospitality and accommodation 
during this visit.  WJM's research is supported by the Canadian government 
through NSERC grant number  OGP0155422.  Additional support was provided by
MITACS, by  CITO through the Centre for Applied Cryptographic Research,
and through NSF-ITR grant number 0112889.

\end{document}